\documentclass[12pt,a4paper]{article}

\usepackage{latexsym,epsfig}
\usepackage{graphicx} %To include graphics
\usepackage[usenames,dvipsnames]{color} %To use color
\usepackage{url} %To use links
\usepackage{soul} %To use highlighting of the text
\usepackage{array} %To use arrays
\usepackage{amssymb}
\usepackage{amsfonts}
\usepackage{amsmath}
\usepackage{hyperref}
\usepackage{verbatim}
\usepackage[sort&compress,square,comma,numbers]{natbib}
\newtheorem{lemma}{Lemma}[section]
\newtheorem{proposition}{Proposition}[section]
\newtheorem{theorem}{Theorem}[section]
\newtheorem{corollary}{Corollary}[section]
\newtheorem{problem}{Problem}

\newtheorem{conjecture}{Conjecture}[section]
\newtheorem{observation}{Observation}[section]
\newcommand{\EndProof}{\hspace{\stretch{1}} $\Box$}

\newcommand{\pr}{\noindent{\bf Proof.}\ }

\newcommand{\G}{\Gamma}

\newcommand{\C}{\mathcal{C}}

\newcommand{\rep}{\mathrm{rep}}

\newcommand{\M}{\mathrm{M}}
\newcommand{\N}{\mathrm{N}}

\title{On large bipartite graphs of diameter 3}
\author{Ramiro Feria-Pur\'on$^{1,}$\footnote{\href{mailto:Ramiro.Feria-Puron@uon.edu.au}{Ramiro.Feria-Puron@uon.edu.au} (Corresponding author)}
\and
Mirka Miller$^{1,2,3,4,}$\footnote{\href{mailto:mirka.miller@newcastle.edu.au}{mirka.miller@newcastle.edu.au}}
\and
Guillermo Pineda-Villavicencio$^{5,}$\footnote{\href{mailto:work@guillermo.com.au}{work@guillermo.com.au}}\and
$\phantom{0}^{1}$\emph{\small School of
Electrical Engineering and Computer Science} \vspace{-1.5mm}\\
\emph{\small The University of Newcastle, Australia}\and
$\phantom{0}^{2}$\emph{\small Department of
Mathematics}\vspace{-1.5mm} \\\emph{\small University of West
Bohemia, Czech Republic}\and
$\phantom{0}^{3}$\emph{\small Department of Computer Science}\vspace{-1.5mm} \\\emph{\small King's College London, UK}\and
$\phantom{0}^{4}$\emph{\small Department of Mathematics}\vspace{-1.5mm} \\\emph{\small ITB Bandung, Indonesia} \and
$\phantom{0}^{5}$\emph{\small Centre for Informatics and Applied Optimization} \vspace{-1.5mm}\\
\emph{\small University of Ballarat, Australia}}

\begin{document}
\maketitle

\begin{abstract}
\noindent
We consider the bipartite version of the {\it degree/diameter problem}, namely, given natural numbers $d\ge2$ and $D\ge2$, find the maximum number $\N^b(d,D)$ of vertices  in a bipartite graph of maximum degree $d$ and diameter $D$. In this context, the bipartite Moore bound $\M^b(d,D)$ represents a general upper bound for $\N^b(d,D)$. Bipartite graphs of order $\M^b(d,D)$ are very rare, and determining $\N^b(d,D)$ still remains an open problem for most $(d,D)$ pairs. 

\noindent
This paper is a follow-up to our earlier paper \cite{FPV12}, where a study on bipartite $(d,D,-4)$-graphs (that is, bipartite graphs of order  $\M^b(d,D)-4$) was carried out. Here we  first present some structural properties of bipartite $(d,3,-4)$-graphs, and later prove there are no bipartite $(7,3,-4)$-graphs. This result implies that the known bipartite $(7,3,-6)$-graph is optimal, and therefore $\N^b(7,3)=80$. Our approach also bears a proof of the uniqueness of the known bipartite $(5,3,-4)$-graph, and the non-existence of bipartite $(6,3,-4)$-graphs.

\noindent
In addition, we discover three new largest known bipartite (and also vertex-transitive) graphs of degree $11$, diameter $3$ and order $190$, result which improves by $4$ vertices the previous lower bound for $\N^b(11,3)$.

\end{abstract}



%										INTRODUCTION

\section{Introduction}

Consider the \emph{degree/diameter problem for bipartite graphs}, stated as follows:

\begin{itemize}
\item[] Given natural numbers $d\ge2$ and $D\ge2$, find the largest possible number $\N^b(d,D)$ of vertices in a bipartite graph of maximum degree $d$ and diameter $D$.
\end{itemize}

It is well known that an upper bound for $\N^b(d,D)$ is given by the {\it bipartite Moore bound} $\M^b(d,D)$, defined below:
\[
\M^b(d,D)=2\left(1+(d-1)+\dots+(d-1)^{D-1}\right).
\]
Bipartite graphs of degree $d$, diameter $D$ and  order $\M^b(d,D)$ are called {\it bipartite Moore graphs}. Bipartite Moore graphs are very scarce; when $d\ge 3$ and $D\ge 3$ they may only exist for $D=3,4$ or $6$ (see \cite{FH64}). It has also turned out to be very difficult to determine $\N^b(d,D)$ even for particular instances; in fact, with the exception of $\N^b(3,5)=\M^b(3,5)-6$ settled in \cite{Jor93},  the other known values of $\N^b(d,D)$ are those for which a bipartite Moore graph is known to exist.

Research in this area falls then into two main directions. On one hand, the efforts to improve the upper bounds for $\N^b(d,D)$ by studying the existence or otherwise of bipartite graphs of  maximum degree $d$, diameter $D$ and order $\M^b(d,D)-\epsilon$ for small $\epsilon>0$ (that is, bipartite $(d,D,-\epsilon)$-graphs, where the parameter $\epsilon$ is called the \emph{defect}). On the other hand, the studies to improve the lower bounds for $\N^b(d,D)$ by constructing ever larger bipartite graphs with given maximum degree and diameter.  In spite of these efforts and the wide range of techniques and approaches used to tackle these problems (see \cite{MS05a}), in most cases there is still a significant gap between the current lower and upper bound for $\N^b(d,D)$.

In this paper we restrict ourselves to the case of bipartite graphs of diameter $3$, and present some modest contributions in both directions. When $D=3$ there is a bipartite Moore graph whenever $d-1$ is a prime power (namely, the incidence graphs of projective planes); however, there is no Moore bipartite graph of diameter $3$ for $d=7$ (\cite{Rys82}) or $d=11$ (\cite{LTS89}). The existence of Moore bipartite graphs of diameter $3$ for other degrees remains an open problem. In \cite{DJMP2} the authors proved that bipartite $(d,3,-2)$-graphs may only exist for certain values of $d$; in particular, they do not exist for $d=7$.

The results and ideas exposed here are, in a great extent, a continuation of the precursory work initiated in \cite{FPV12}. We provide structural properties for bipartite $(d,3,-4)$-graphs and, most important,  prove the non-existence of bipartite $(7,3,-4)$-graphs. Such outcome implies that the only known bipartite $(7,3,-6)$-graph -- found by Paul Hafner and independently by Eyal Loz (\cite{LH12}) -- is optimal, and therefore $\N^b(7,3)=80$. This is just the second  value settled for $\N^b(d,D)$ other than a bipartite Moore bound.
Our approach can also be used to show the uniqueness of the known bipartite $(5,3,-4)$-graph, as well as the non-existence of bipartite $(6,3,-4)$-graphs.

Finally, we also find three largest known bipartite (and vertex-transitive) graphs of degree $11$ and diameter $3$. This settles $190\le \N^b(11,3)$, which improves by $4$ vertices the previous lower bound for $\N^b(11,3)$. Adjacency lists of these graphs are available at \url{http://guillermo.com.au/wiki/List_of_Publications} under the name of this paper.

We conclude this introduction by depicting all the known bipartite $(d,3,-4)$ graphs. Figure \ref{fig:Bipartite(3,3,-4)} shows all the bipartite $(3,3,-4)$-graphs, Figure \ref{fig:Bipartite(4,3,-4)} all the bipartite $(3,3,-4)$-graphs, and Figure \ref{fig:Bipartite(5,3,-4)} the -- after this paper unique -- bipartite $(5,3,-4)$-graph.

\begin{figure}[phtb]
\begin{center}
\makebox[\textwidth][c]{\includegraphics[scale=.75]{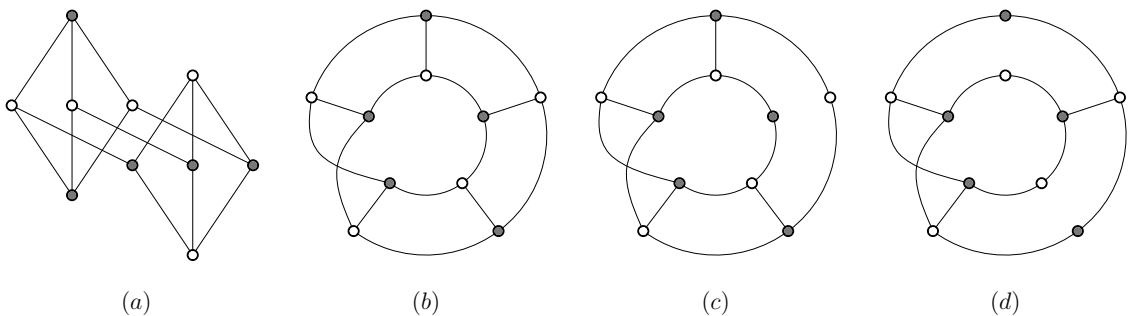}}
\caption{All the bipartite $(3,3,-4)$-graphs.}
\label{fig:Bipartite(3,3,-4)}
\end{center}
\end{figure}

\begin{figure}[phtb]
\begin{center}
\makebox[\textwidth][c]{\includegraphics[scale=.75]{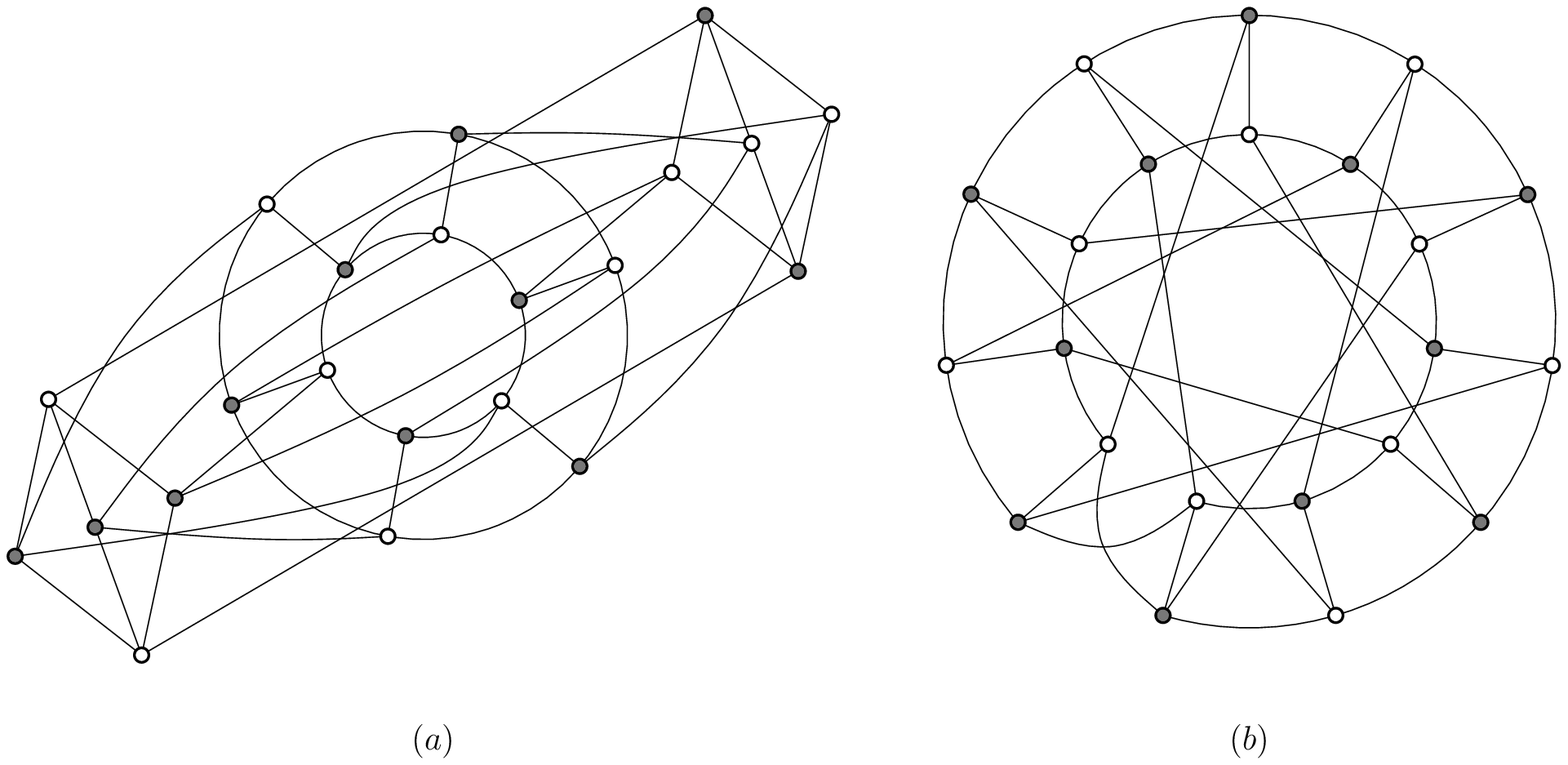}}
\caption{All the bipartite $(4,3,-4)$-graphs.}
\label{fig:Bipartite(4,3,-4)}
\end{center}
\end{figure}

\begin{figure}[phtb]
\begin{center}
\makebox[\textwidth][c]{\includegraphics[scale=.75]{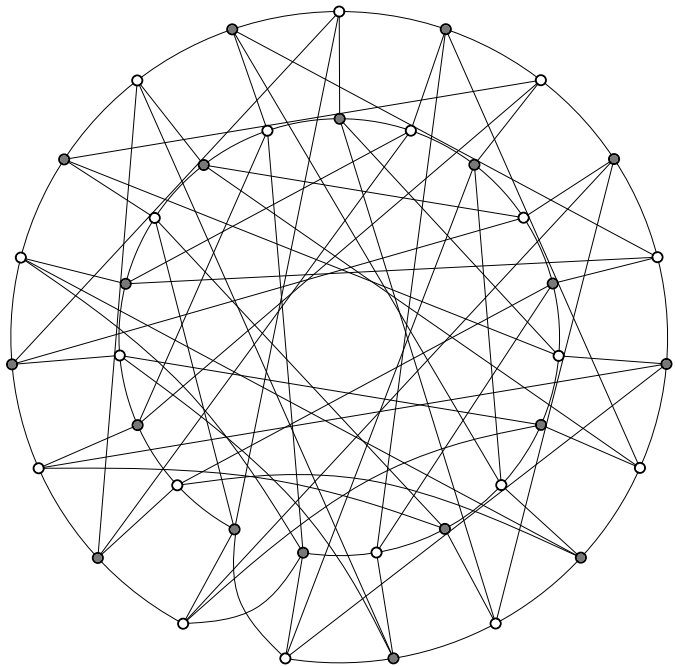}}
\caption{The unique bipartite $(5,3,-4)$-graph.}
\label{fig:Bipartite(5,3,-4)}
\end{center}
\end{figure}

%										NOTATION

\section{Notation and Terminology}

Our notation and terminology follows from \cite{FPV12}, which is standard and consistent with that used in \cite{Die05}.

All graphs considered are simple. The vertex set of a graph $\G$ is denoted by $V(\G)$, and its edge set by $E(\G)$. For an edge $e=\{x,y\}$ we write $x\sim y$. The set of edges in a graph $\G$ joining a vertex $x$ in $X\subseteq V(\G)$ to a vertex $y$ in $Y\subseteq V(\G)$ is denoted by $E(X,Y)$. A vertex of degree at least 3 is called a \emph{branch vertex} of $\G$.

A cycle of length $k$ is called a {\it $k$-cycle}. In a bipartite $(d,D,-4)$-graph we call a cycle of length at most $2D-2$ a {\it short cycle}. If two short cycles $C^1$ and $C^2$  are non-disjoint we say that $C^1$ and $C^2$ are \emph{neighbors}.

For a vertex $x$ lying on a short cycle $C$, we denote by $\rep^C(x)$ the vertex $x'$ in $C$ such that $d(x,x')=D-1$, where $d(x,x')$ denotes the distance between $x$ and $x'$. In this case, we say $x'$ is a \emph{repeat} of $x$ in $C$ and vice versa, or simply that $x$ and $x'$ are \emph{repeats} in $C$. A \emph{closed set of repeats} in a bipartite $(d,D,-4)$-graph $\G$ is a subset of $V(\G)$ which is closed under the repeat relation. A closed set of repeats is \emph{minimal} if it does not have a proper closed subset of repeats.

Finally, we introduce some special graphs. The union of three independent paths of length $t$ with common endvertices is denoted by $\Theta_t$. For an integer $m \ge 5$,  $\Phi_m$ denotes the bipartite graph with vertex set $V=\{x_i | 0\le i \le m-1\} \cup \{y_i | 0\le i \le m-1\}$ and edge set $E=\{x_i\sim y_i, x_i\sim y_{i+1}, x_i\sim y_{i-1} | 0\le i \le m-1\}$. Note that $\Phi_m$ is vertex-transitive. Throughout this paper we do addition modulo $m$ on the vertex subscripts of a $\Phi_m$.

%										PRELIMINARIES

\section{Preliminaries}
\label{sec:Preliminaries}

We begin with the regularity condition for bipartite graphs with small defect.

\begin{proposition}[\cite{DJMP2}]
\label{prop:regularity1} For $\epsilon<1+(d-1)+(d-1)^2+\ldots+(d-1)^{D-2}$, $d\ge 3$ and $D\ge 3$, a bipartite $(d,D,-\epsilon)$-graph is regular.
\end{proposition}

\begin{proposition}[\cite{DJMP2}]
\label{prop:regularity2} For $\epsilon<2\left((d-1)+(d-1)^3+\ldots+(d-1)^{D-2}\right)$,$d\ge 3$ and odd $D\ge 3$, a bipartite $(d,D,-\epsilon)$-graph is regular.
\end{proposition}

In particular, we will implicitly use the fact that a bipartite $(d,3,-4)$-graph with $d\ge4$ must be regular, and therefore its partite sets must have the same cardinality. Also note that, since bipartite $(d,3,-\epsilon)$ graphs with $d \ge 4$ and $\epsilon =3,5$ are not regular, the above propositions imply their non-existence.

From the paper \cite{FPV12} we borrow the following results:

%Proposition $4.1$, the Saturating Lemma, the Repeat Cycle Lemma, and the following propositions:

\begin{proposition}[\cite{FPV12}]
\label{prop:GirthBipartite(d,D,-4)} The girth of a regular bipartite $(d,D,-4)$-graph $\G$ with $d\ge3$ and $D\ge 3$ is $2D-2$. Furthermore, any vertex $x$ of $\G$ lies on the short cycles specified below and no other short cycle, and we have the following cases:

\noindent {\bf $x$ is contained in {\bf exactly three $(2D-2)$-cycles}}. Then
\begin{itemize}
\item[{\rm ($i$)}] $x$ is a branch vertex of one $\Theta_{D-1}$, or
\end{itemize}
{\bf $x$ is contained in {\bf two $(2D-2)$-cycles}}. Then
\begin{itemize}
\item[{\rm ($ii$)}] $x$ lies on {\bf exactly two $(2D-2)$-cycles}, whose
intersection is a $\ell$-path with $\ell\in \{0,\ldots,D-1\}$.
\end{itemize}
\end{proposition}

As in \cite{FPV12}, often our arguments revolve around the identification of the elements in the set $S_x$ of short cycles containing a given vertex $x$; we call this process \emph{saturating} the vertex $x$. A vertex $x$ is called \emph{saturated} if the elements in $S_x$ have been completely identified.

\begin{lemma}[\cite{FPV12}, Saturating Lemma]
\label{lemm:SaturatingLemma} Let $\C$ be a $(2D-2)$-cycle in a bipartite $(d,D,-4)$-graph $\G$ with $d \ge 4$ and $D \ge 3$, and $\alpha,\alpha'$ two vertices in $\C$ such that $\alpha'=\rep^{\C}(\alpha)$. Let $\gamma$ be a neighbor of $\alpha$ not contained in $\C$, and $\mu_1,\mu_2,\ldots,\mu_{d-2}$ the neighbors of $\alpha'$ not contained in $\C$. Suppose there is no short cycle in $\G$ containing the edge $\alpha\sim\gamma$ and intersecting $\C$ at a path of length greater than $D-3$.

Then, in $\G$ there exist  a vertex $\mu\in\{\mu_1,\mu_2,\ldots,\mu_{d-2}\}$ and a short cycle $\C^1$ such that $\gamma$ and $\mu$ are repeats in $\C^1$, and $\C\cap\C^1=\emptyset$.
\end{lemma}

\begin{lemma}[\cite{FPV12}, Repeat Cycle Lemma]
\label{lemm:RepeatCycleLemma}
Let $C$ be a short cycle in a bipartite $(d,D,-4)$-graph $\G$ with $d \ge 4$ and $D \ge 3$, $\{C^1,C^2,\ldots ,C^k\}$ the set of neighbors of $C$, and  $I_i=C^i\cap C$ for $1\le i\le k$. Suppose at least one $I_j$, for $j\in\{1,\ldots, k\}$, is a path of length smaller than $D-2$.  Then there is an additional short cycle $C'$ in $\G$ intersecting $C^i$ at $I'_i=\rep^{C^i}(I_i)$, where $1\le i\le k$.
\end{lemma}

\begin{proposition}[\cite{FPV12}]
\label{prop:CyclePartition}
The set $S(\G)$ of short cycles in a bipartite $(d,D,-4)$-graph $\G$ with $d \ge 3$ and $D \ge 3$ can be partitioned into sets $S_{D-1}(\G)$, $S_{D-2}(\G)$ and $S_{D-3}(\G)$, where

\begin{itemize}
\item[] $S_{D-1}(\G)$ is the set of short cycles in $\G$ whose intersections with neighbor cycles are $(D-1)$-paths,
\item[] $S_{D-2}(\G)$ is the set of short cycles in $\G$ whose intersections with neighbor cycles are $(D-2)$-paths, and
\item[] $S_{D-3}(\G)$ is the set of short cycles in $\G$ whose intersections with neighbor cycles are paths of length at most $D-3$.
\end{itemize}
\end{proposition}

\begin{proposition}[\cite{FPV12}]
\label{prop:VertexPartition}
The set $V(\G)$ of vertices in a bipartite $(d,D,-4)$-graph $\G$ with $d \ge 4$ and $D \ge 3$ can be partitioned into sets $V_{D-1}(\G)$, $V_{D-2}(\G)$ and $V_{D-3}(\G)$, where

\begin{itemize}
\item[] $V_{D-1}(\G)$ is the set of vertices contained in cycles of  $S_{D-1}(\G)$,
\item[] $V_{D-2}(\G)$ is the set of vertices contained in cycles of  $S_{D-2}(\G)$,
\item[] $V_{D-3}(\G)$ is the set of vertices contained in cycles of  $S_{D-3}(\G)$,
\end{itemize}

and $S_{D-1}(\G)$, $S_{D-2}(\G)$, $S_{D-3}(\G)$ are defined as in Proposition \ref{prop:CyclePartition}.
\end{proposition}

%										ON BIPARTITE GRAPHS OF D = 3 AND DEFECT 4 (Observations)

\subsection{On bipartite graphs of diameter $3$ and defect $4$}

In this section we present additional structural properties for bipartite graphs of diameter 3 and defect $4$.

Let $\G$ be a bipartite $(d,3,-4)$-graphs with $d\ge4$. We set $\G_i=\bigcup_{C \in S_i(\G)} C$  for $i=0,1,2$. Note that $\G_2$ is the union of all graphs in $\G$ isomorphic to $\Theta_2$; these graphs are pairwise disjoint, so they are the connected components of $\G_2$. In addition, $\G_1$ is the union of all graphs in $\G$ isomorphic to some $\Phi_m$ for $m\ge5$; similarly, these $\Phi_m$ are the connected components of $\G_1$. 

If $G$ is a connected component in $\G_2 \cup \G_1 \cup \G_0$ then $V(G)$ is a closed set of repeats.  The branch vertices of a $\Theta_2 \subset \G_2$ constitute a minimal closed set of repeats, as well as its non-branch vertices. In the case of a $\Phi_m \subset \G_1$, the vertices $x_i$'s form a minimal closed set of repeats, the same as the vertices $y_i$'s. According to the Repeat Cycle Lemma, every minimal closed set of repeats in $\G_0$ contains exactly $4$ vertices. Observe that all vertices in a minimal closed set of repeats in $\G$ belong to the same partite set.

Some further observations about $\G$ follow from the systematic application of the Saturating Lemma:

\begin{observation}
\label{obs:noedgeG2G2}
Let $\G$ be a bipartite $(d,3,-4)$-graph with $d\ge4$. There is no edge in $\G$ joining a branch vertex in $\G_2$ to a non-branch vertex of a different connected component of $\G_2$.
\end{observation}

\pr

\begin{figure}[!ht]
\begin{center}
\makebox[\textwidth][c]{\includegraphics[scale=.8]{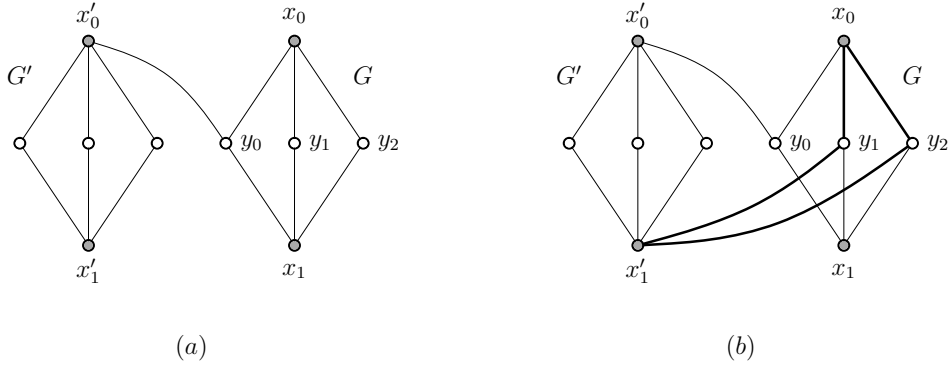}}
\caption{Auxiliary figure for Observation \ref{obs:noedgeG2G2}}
\label{fig:noedgeG2G2}
\end{center}
\end{figure}

Let $G,G'$ be two connected components in $\G_0$ such that a branch vertex $x'_0$ in $G'$ is adjacent to a non-branch vertex $y_0$ in $G$. Let $x'_1,x_0,x_1,y_1,y_2$ be as in Figure \ref{fig:noedgeG2G2} $(a)$. We apply the Saturating Lemma (by mapping the cycle $x_0y_0x_1y_1x_0$ to $\C$, $y_0$ to $\alpha$, $y_1$ to $\alpha'$ and $x'_0$ to $\gamma$), and obtain that $y_1$ is adjacent to $x'_1$. Similarly, $y_2$ is also adjacent to $x'_1$ (see Figure \ref{fig:noedgeG2G2} $(b)$), but then there is a fourth short cycle $x_0y_1x'_1y_2x_0$ in $\G$ containing $x_0$, a contradiction. \EndProof

\begin{observation}
\label{obs:noedgeG2G1}
Let $\G$ be a bipartite $(d,3,-4)$-graph with $d \ge 4$. There is no edge in $\G$ joining a branch vertex in $\G_2$ to a vertex in $\G_1$.
\end{observation}

\pr

\begin{figure}[!ht]
\begin{center}
\makebox[\textwidth][c]{\includegraphics[scale=.8]{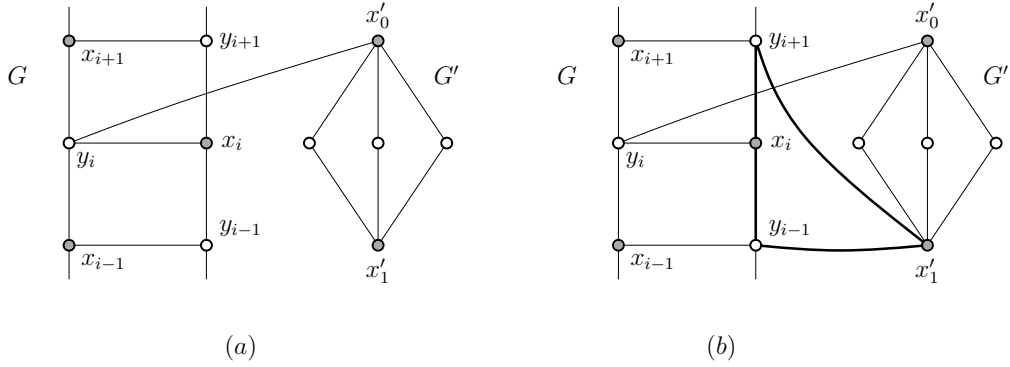}}
\caption{Auxiliary figure for Observation \ref{obs:noedgeG2G1}}
\label{fig:noedgeG2G1}
\end{center}
\end{figure}

Let $G,G'$ be two connected components of $\G_1$ and $\G_2$ respectively, such that a branch vertex $x'_0$ in $G'$ is adjacent to a  vertex $y_i$ in $G$. Let $x'_1, y_{i+1}, y_{i-1}, x_i, x_{i+1}, x_{i-1}$ be as in Figure \ref{fig:noedgeG2G1} $(a)$. We apply the Saturating Lemma (by mapping cycle $y_ix_{i-1}y_{i-1}x_iy_i$ to $\C$, $y_i$ to $\alpha$, $y_{i-1}$ to $\alpha'$ and $x'_0$ to $\gamma$), and obtain that $y_{i-1}$ is adjacent to $x'_1$. Similarly, $y_{i+1}$ is also adjacent to $x'_1$ (see Figure \ref{fig:noedgeG2G1} $(b)$). But then, there is a third short cycle $y_{i+1}x_iy_{i-1}x'_1y_{i+1}$ in $\G$ containing $x_i$, a contradiction. \EndProof

\begin{observation}
\label{obs:noedgeG2G0}
Let $\G$ be a bipartite $(d,3,-4)$-graph with $d \ge 4$. There is no edge in $\G$ joining a non-branch vertex in $\G_2$ to a vertex in $\G_0$.
\end{observation}

\pr

\begin{figure}[!ht]
\begin{center}
\makebox[\textwidth][c]{\includegraphics[scale=.8]{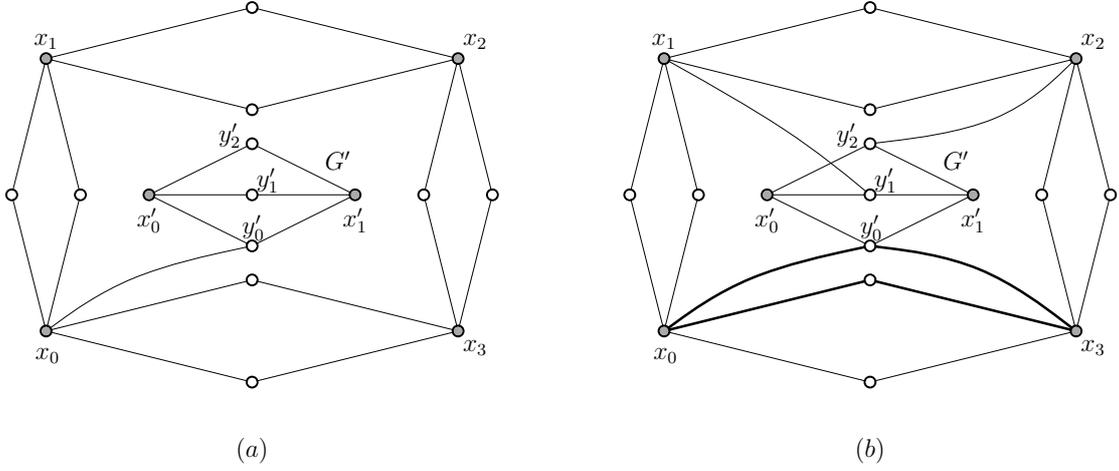}}
\caption{Auxiliary figure for Observation \ref{obs:noedgeG2G0}}
\label{fig:noedgeG2G0}
\end{center}
\end{figure}

Let $G'$ be a connected component in $\G_2$ with a non-branch vertex $y'_0$ adjacent to a vertex $x_0$ in $\G_0$. Let $\{x_0,x_1,x_2,x_3\}$ be the minimal closed set of repeats containing $x_0$ ($x_2$ not being a repeat of $x_0$), and let the vertices $x'_0,x'_1,y'_1,y'_2$ be as in Figure \ref{fig:noedgeG2G0} $(a)$. We first apply the Saturating Lemma (by mapping the cycle $x'_0y'_0x'_1y'_1x'_0$ to $\C$, $y'_0$ to $\alpha$, $y'_1$ to $\alpha'$ and $x_0$ to $\gamma$), and obtain that $y'_1$ is adjacent to a repeat of $x_0$ (say $x_1$). Similarly, mapping the cycle $x'_0y'_0x'_1y'_1x'_0$ to $\C$, $y'_1$ to $\alpha$, $y'_2$ to $\alpha'$ and $x_1$ to $\gamma$, we obtain that $y'_2$ is adjacent to $x_2$ (as it cannot be adjacent to $x_0$). Analogously,  $y'_0$ is adjacent to $x_3$ (see Figure \ref{fig:noedgeG2G0} $(b)$), but then there is a third short cycle in $\G$ containing $x_0$, a contradiction.\EndProof

\begin{observation}
\label{obs:G1}
Let $\G$ be a bipartite $(d,3,-4)$-graph with $d \ge 4$, and $G=\Phi_m$ a connected component in $\G_1$. Given $x_i\in G$, if $x_i\sim y_j \in E(\G)$ for some $j$ then $x_{i+k}\sim y_{j+k} \in E(\G)$ for every $k$.
\end{observation}

\pr

This clearly holds when $j \in \{i,i+1,i-1\}$; see the description of $\Phi_m$.

Suppose $j \not \in \{i,i+1,i-1\}$. Since all the vertices in $G$ are saturated, we have $|i-j|\ge4$. According to the Saturating Lemma (by mapping the cycle $x_iy_ix_{i+1}y_{i+1}x_i$ to $\C$, $x_i$ to $\alpha$, $x_{i+1}$ to $\alpha'$ and $y_j$ to $\gamma$) we have either $x_{i+1}\sim y_{j+1} \in E(\G)$ or $x_{i+1}\sim y_{j-1} \in E(\G)$. But in case $x_{i+1}\sim y_{j-1} \in E(\G)$, it is easy to see that, by repeatedly applying the Saturating Lemma (to the cycles $x_{i+p}y_{i+p}x_{i+p+1}y_{i+p+1}x_{i+p}$ for $p=1,2,\ldots$) we obtain there is an edge $x_r\sim y_s$ in $\G$ such that $2\le|r-s|\le3$, which is not possible. Thus $x_{i+1}\sim y_{j+1} \in E(\G)$ and, by induction, $x_{i+k}\sim y_{j+k} \in E(\G)$ for every $k$.\EndProof

\begin{observation}
\label{obs:G1G1}
Let $\G$ be a bipartite $(d,3,-4)$-graph with $d\ge3$, and $G,G'$ two connected components in $\G_1$ of order $2m$ and $2m'$ respectively ($m \le m'$). Suppose there is at least one edge in $\G$ joining a vertex in $G$ to a vertex in $G'$. Then $m'=km$, with $1 \le k \le d-3$.
\end{observation}

\pr 

Denote the vertices of $G=\Phi_m$ by $x_0,\ldots,x_{m-1},y_0,\ldots,y_{m-1}$, and the vertices of $G'=\Phi_{m'}$ by $x'_0,\ldots,x'_{m'-1},y'_0,\ldots,y'_{m'-1}$. With an appropriate labelling we may assume there is an edge $x_0\sim y'_0$ in $\G$ and, by the Saturating Lemma (on the cycle $x_0y_0x_1y_1x_0$), also an edge $x_1\sim y'_1$ in $\G$. 

Suppose $m'=km+r$, with $1 \le r \le m-1$ and $k \ge 1$. Then, by repeatedly applying the Saturating Lemma on the cycles $x_iy_ix_{i+1}y_{i+1}x_i$ with $i=1,\ldots,m-1$, we find the edges $x_i\sim y'_i$ for $i=2,\ldots,m$ are all present in $\G$. In particular, $y'_m$ is a neighbor of $x_0$ and, inductively, the vertices $y'_{2m},\ldots, y'_{km}, y'_{m-r}, y'_{2m-r},\ldots$ also are. But similarly, $x_{m-r}$ has also neighbors  $y'_{m-r}$ and $y'_{2m-r}$; this way, we obtain there is in $\G$ a third short cycle $x_0y'_{m-r}x_{m-r}y'_{2m-r}x_0$ containing $x_0$, a contradiction.

Since a vertex in $G$  has at least $3$ neighbors in $G$, it follows that $k \le d-3$.
 \EndProof

\begin{observation}
\label{obs:cardinalG0}
Let $\G$ be a bipartite $(7,3,-4)$-graph. If $\G_0 \neq \emptyset$ then $|\G_0|=8k$, with $k \ge 3$.
\end{observation}

\pr

If $t$ is the number of short cycles in $\G_0$ then, by a simple counting argument, $\G_0$ has $2t$ vertices, half of them in each partite set. Recall that $V(\G_0)$ is a closed set of repeats. Since a minimal closed set of repeats in $\G_0$ contains exactly $4$ vertices belonging to the same partite set, we have $t=4k$ and then $|\G_0|=8k$.

\begin{figure}[!ht]
\begin{center}
\makebox[\textwidth][c]{\includegraphics[scale=.8]{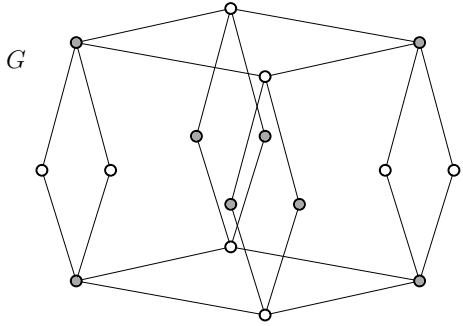}}
\caption{Auxiliary figure for Observation \ref{obs:cardinalG0}}
\label{fig:cardinalG0}
\end{center}
\end{figure}

Also, the Repeat Cycle Lemma ensures that the graph $G$ depicted in Figure \ref{fig:cardinalG0} is a subgraph of $\G_0$. Since any vertex in $\G_0$ must have at least $4$ neighbors in $\G_0$, we have $|\G_0| > 16$ and $k\ge3$.
\EndProof

\begin{observation}
\label{obs:edgeG2G1}
Let $\G$ be a bipartite $(d,3,-4)$-graph with $d\ge4$, $G$ a connected component in $\G_2$, and $G'$ a connected component in $\G_1$ of order $2m'$. Suppose there is in $\G$ at least one edge joining a vertex in $G$ to a vertex in $G'$. Then $m'=3k$ with $2 \le k \le d-2$.
\end{observation}

\pr

Let $x_0,x_1,x_2$ be the non-branch vertices of $G=\Theta_2$, and denote by $x'_0,\ldots,x'_{m'-1},y'_0,\ldots,y'_{m'-1}$ the vertices of $G'=\Phi_{m'}$.

By Observation \ref{obs:noedgeG2G1} any edge between $\G_2$ and $\G_1$ involves only non-branch vertices of $\G_2$. We may assume there are edges $x_0\sim y'_0$ and $x_1\sim y'_1$ in $\G$. Suppose $m'=3k+r$, with $1 \le r \le 2$ and $k \ge 1$. Then, by repeatedly applying the Saturating Lemma on the three short cycles of $G$, we obtain that $x_0$ has neighbors $y'_0, y'_3, y'_6,\ldots, y'_{3k}, y'_{3-r}, y'_{6-r},\ldots$ But similarly, $x_{3-r}$ has also neighbors $y'_{3-r}$ and $y'_{6-r}$; hence, we obtain there is in $\G$ a third short cycle $x_0y'_{3-r}x_{3-r}y'_{6-r}x_0$ containing $x_0$, a contradiction.

Since each $x_i$ has $2$ neighbors in $G$ and $m' \ge 5$, it follows that $2 \le k \le d-2$.
 \EndProof

\begin{observation}
\label{obs:edgeG0G1}
Let $\G$ be a bipartite $(d,3,-4)$-graph with $d\ge4$ and $G'$ a connected component in $\G_1$ of order $2m'$. Suppose there is in $\G$ an edge joining a vertex in $\G_0$ to a vertex in $G'$. Then $m'=4k$ with $2 \le k \le d-4$.
\end{observation}

\pr

Let $x_0\in V(\G_0)$ and let $\{x_0,x_1,x_2,x_3\}$ be the minimal closed set of repeats containing $x_0$ ($x_2$ not being a repeat of $x_0$). Denote by $x'_0,\ldots,x'_{m'-1},y'_0,\ldots,y'_{m'-1}$ the vertices of $G'=\Phi_{m'}$. We may assume there are edges $x_0\sim y'_0$ and $x_1\sim y'_1$ in $\G$. Suppose $m'=4k+r$, with $1 \le r \le 3$ and $k \ge 1$. Then, by repeatedly applying the Saturating Lemma on the cycles $x'_iy'_ix'_{i+1}y'_{i+1}x'_i$ ($i=1,2,\ldots$) of $G'$, we obtain that $x_0$ has neighbors $y'_0, y'_4, y'_8,\ldots, y'_{4k}, y'_{4-r}, y'_{8-r},\ldots$ But analogously, $x_{4-r}$ has also neighbors $y'_{4-r}$ and $y'_{8-r}$; hence, we obtain there is in $\G$ a third short cycle $x_0y'_{4-r}x_{4-r}y'_{8-r}x_0$ containing $x_0$, a contradiction.

Since $x_0$ has at least $4$ neighbors in $\G_0$ and $m' \ge 5$, it follows that $2 \le k \le d-4$.
 \EndProof

The statements in Observations \ref{obs:noedgeG2G2}, \ref{obs:noedgeG2G0}, \ref{obs:G1G1}, \ref{obs:edgeG2G1} and \ref{obs:edgeG0G1} are better summarized in the following, more compact assertion.

\begin{proposition}
\label{prop:edgeMM'}
Let $\G$ be a bipartite $(d,3,-4)$-graph with $d\ge4$, and $M,M'$ two minimal closed set of repeats in $\G$ such that $E(M,M')\neq \emptyset$. Then $|M|$ divides $|M'|$ or $|M'|$ divides $|M|$, except when $M \cup M'$ is the set of five the vertices in a $\Theta_2$ .
\end{proposition}\EndProof

%										Nb(7,3)=80

\section{Non-existence of bipartite $(7,3,-4)$-graphs}
\label{sec:Nb(7,3)=80}

In this section we prove that there are no bipartite $(7,3,-4)$-graphs, and consequently that $\N^b(7,3)=80$.

\begin{proposition}
\label{prop:G2}
Let $\G$ be a bipartite $(7,3,-4)$-graph. Then $\G_2$ cannot be a spanning subgraph of $\G$.
\end{proposition}

\pr

Since the connected components of $\G_2$ are graphs isomorphic to $\Theta_2$, we have that $5$ must divide $|\G_2|=82$, a contradiction.\EndProof

\begin{proposition}
\label{prop:G1}
Let $\G$ be a bipartite $(7,3,-4)$-graph. Then $\G_1$ cannot be a spanning subgraph of $\G$.
\end{proposition}

\pr

This is a computer-assisted proof.

Suppose that $\G_1$ contains exactly one connected component $G$, isomorphic to $\Phi_{41}$. Denote by $x_0,\ldots,x_{40},y_0,\ldots,y_{40}$ the vertices of $G$. By virtue of Observation \ref{obs:G1}, if the vertex $x_0$ has neighbors $y_0,y_1,y_{-1},y_{i_1},y_{i_2},y_{i_3},y_{i_4}$ in $G$ then $x_k$ has neighbors $y_k,y_{k+1},y_{k-1},y_{k+i_1},y_{k+i_2},y_{k+i_3},y_{k+i_4}$ for every $k$. Exhaustive computer search through the feasible choices for the vertices $y_{i_1},y_{i_2},y_{i_3},y_{i_4}$ yields no graph of diameter $3$, and so there is more than one conected component in $\G_1$.

Now suppose then that $\G_1$ has exactly $n$ connected components $G_1,G_2,\ldots,G_n$, isomorphic to $\Phi_{m_1},\Phi_{m_2},\ldots,\Phi_{m_n}$, respectively. Note that $5\le m_i \le 36$, $2\le n\le 8$ and $m_1+\ldots+m_n=41$. We define the graph $H(G_1,G_1,\ldots,G_n)$ as follows: every $G_i$ contracts to a vertex $v_i$ in $H$, and there is an edge $v_i-v_j$  in $H$ if and only if -- according to Observation \ref{obs:G1G1} -- there could be an edge from $G_i$ to $G_j$ in $\G$ (that is, if $m_i$ divides $m_j$ or vice versa). Clearly, if $\G$ has diameter $3$ then $H$ has diameter at most $2$. However, we could verify that none of the feasible values for $n$ and the $m_i$'s yields a graph $H$ of diameter at most $2$.

Consequently, $V(\G_1)$ cannot span $\G$. \EndProof

\begin{proposition}
\label{prop:G0}
Let $\G$ be a bipartite $(7,3,-4)$-graph. Then $\G_0$ cannot be a spanning subgraph of $\G$.
\end{proposition}

\pr

From Observation \ref{obs:cardinalG0} we have $82=|\G_0|=8k$, a contradiction. \EndProof

\begin{proposition}
\label{prop:G2G1}
Let $\G$ be a bipartite $(7,3,-4)$-graph. Then $\G_2 \cup \G_1$ cannot be a spanning subgraph of $\G$.
\end{proposition}

\pr

Suppose $\G_2 \neq \emptyset$ and $\G_1 \neq \emptyset$. On one hand, from a branch vertex in $\G_2$ it is possible to reach in exactly two steps at most $15$ vertices of $\G_1$ (see Observations \ref{obs:noedgeG2G2} and \ref{obs:noedgeG2G1}). Therefore, $|\G_1| \le 30$. On the other hand, from a vertex in $\G_1$ it is possible to reach in exactly two steps at most $8$ branch vertices of $\G_2$, and $|\G_2| \le 40$. This means $|\G| \le 70$, a contradiction. \EndProof

\begin{proposition}
\label{prop:G2G0}
Let $\G$ be a bipartite $(7,3,-4)$-graph. Then $\G_2 \cup \G_0$ cannot be a spanning subgraph of $\G$.
\end{proposition}

\pr

Suppose $\G_2 \neq \emptyset$ and $\G_0 \neq \emptyset$. From a non-branch vertex in $\G_2$ it is possible to reach in two steps at most $8$ vertices of $\G_0$ (see Observations \ref{obs:noedgeG2G2} and \ref{obs:noedgeG2G0}). Therefore, $|\G_0| \le 16$, which contradicts Observation \ref{obs:cardinalG0}. \EndProof

\begin{proposition}
\label{prop:G1G0}
Let $\G$ be a bipartite $(7,3,-4)$-graph. Then $\G_1 \cup \G_0$ cannot be a spanning subgraph of $\G$.
\end{proposition}

\pr

Let $G=\Phi_m$ be a connected component in $\G_1$. We prove that $m$ is even. If $G$ has a neighbor in $\G_0$ then, by Observation \ref{obs:edgeG0G1}, we have $m \in \{8,12\}$. If instead $G$ has no neighbor in $\G_0$ and $m$ is odd, then there must be a connected component $G'$ in $\G_1$ isomorphic to some $\Phi_{m'}$ such that $G$ has a neighbor in $G'$ and $G'$ has a neighbor in $\G_0$. But again we have $m' \in \{8,12\}$ and, according to Observation \ref{obs:G1G1}, $m\ge5$ must be an odd divisor of $m'$, which is not possible.

From the above and Observation \ref{obs:cardinalG0} it follows that $|\G| \equiv 0 \pmod 4$, which contradicts $|\G|=82$.
\EndProof

\begin{proposition}
\label{prop:G2G1G0}
Let $\G$ be a bipartite $(7,3,-4)$-graph. Then $\G_2 \cup \G_1 \cup \G_0$ cannot be a spanning subgraph of $\G$.
\end{proposition}

\pr

Let $\G_i \neq \emptyset$ for $i=0,1,2$.

\noindent {\bf Claim 1.} Every connected component of $\G_1$ has a neighbor in $\G_0$.

\noindent {\bf Proof of Claim 1.}

\begin{figure}[!ht]
\begin{center}
\makebox[\textwidth][c]{\includegraphics[scale=.9]{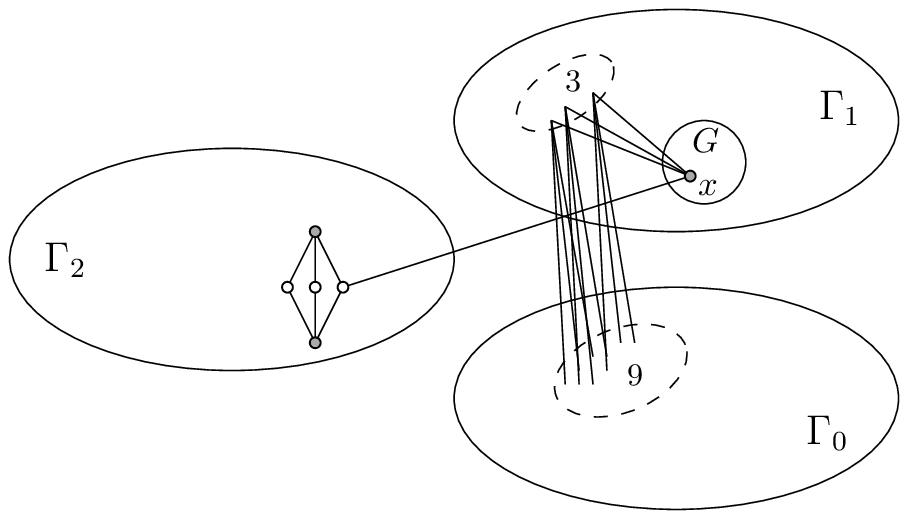}}
\caption{Auxiliary figure for Claim1}
\label{fig:Claim1G2G1G0}
\end{center}
\end{figure}

Suppose there is a connected component $G$ of $\G_1$ with no neighbors in $\G_0$, and take a vertex $x$ in $G$. According to Observation \ref{obs:noedgeG2G1}, $x$ must have at least one non-branch neighbor in $\G_2$ for it can reach in $2$ steps the branch vertices of $\G_2$ belonging to its partite set.  But then from $x$ it is possible to reach at most $9$ vertices of $\G_0$ in exactly $2$ steps (see Figure \ref{fig:Claim1G2G1G0}). This implies $|\G_0|\le18$, which contradicts Observation \ref{obs:cardinalG0}.
\EndProof

\noindent {\bf Claim 2.} Every connected component of $\G_1$ has a neighbor in $\G_2$.

\noindent {\bf Proof of Claim 2.}

\begin{figure}[!ht]
\begin{center}
\makebox[\textwidth][c]{\includegraphics[scale=.9]{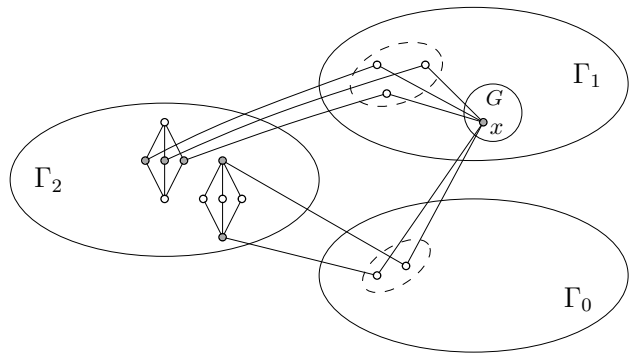}}
\caption{Auxiliary figure for Claim 2}
\label{fig:Claim2G2G1G0}
\end{center}
\end{figure}

Suppose there is a connected component $G$ of $\G_1$ with no neighbors in $\G_2$. First note that $\G_2$ must have the same number of vertices in each partite set of $\G$, so $|\G_2|\ge10$. From a vertex $x$ in $G$ we must reach in exactly two steps at least three non-branch vertices in a connected component of $\G_2$, and other two branch vertices in a different connected component of $\G_2$. However, it is only possible to reach from $x$ at most $4$ of such $5$ vertices (see Figure \ref{fig:Claim2G2G1G0}).
\EndProof

From Claim 1 and Observation \ref{obs:edgeG0G1} we can deduce that if $G=\Phi_m$ is a connected component of $\G_1$ then $m \in \{8,12\}$. But  from Claim 2 and Observation \ref{obs:edgeG2G1} it follows that $m \equiv 0 \pmod{3}$, and therefore  $m=12$. In other words, every connected component of $\G_1$ has $24$ vertices.

In addition, since $|\G_0|\ge24$ and $|\G_1|\ge24$ we have that $|\G_2|\le34$. But $5$ (and hence $10$) must divide $|\G_2|$, and $|\G_0| \equiv |\G_1| \equiv 0 \pmod 8$; consequently, $|\G_2|=10$. 

\begin{figure}[!ht]
\begin{center}
\makebox[\textwidth][c]{\includegraphics[scale=.9]{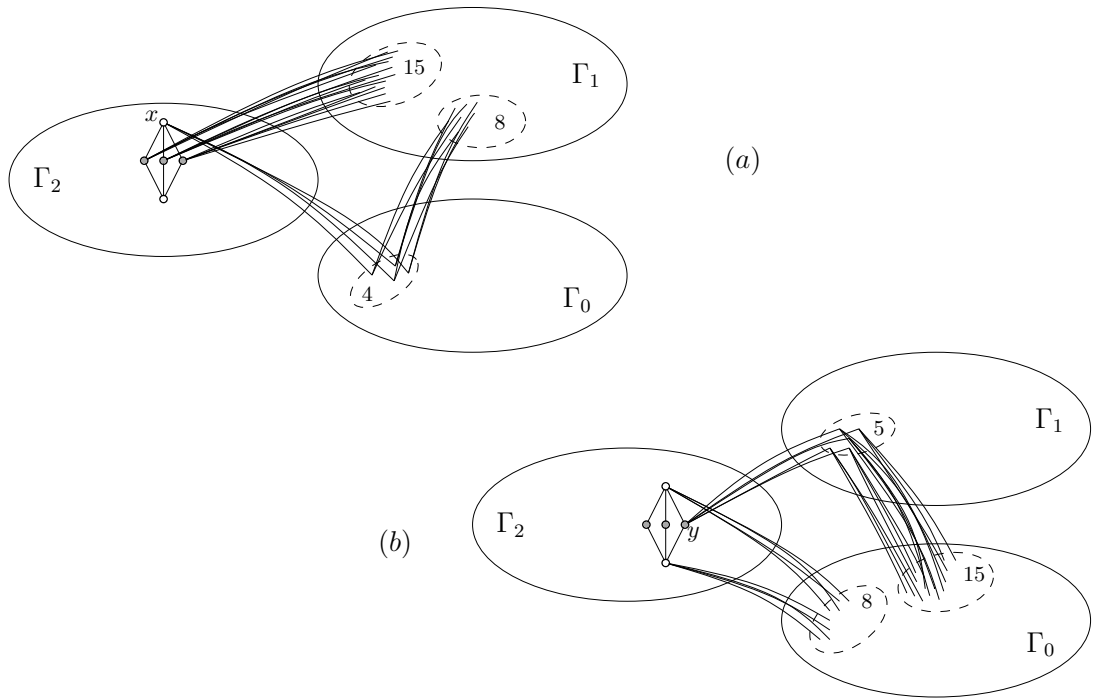}}
\caption{Auxiliary figure for Proposition \ref{prop:G2G1G0}}
\label{fig:G2G1G0}
\end{center}
\end{figure}

To complete the proof we only need to consider two possibilities left. If $|\G_2|=10$, $|\G_1|=48$ and $|\G_0|=24$ then from a branch vertex $x$ in $\G_2$ it is possible to reach in exactly two steps at most $23$ vertices of $\G_1$ in the same partite set as $x$, a contradiction (see Figure \ref{fig:G2G1G0} $(a)$). Similarly, if $|\G_2|=10$, $|\G_1|=24$ and $|\G_0|=48$ then from a non-branch vertex $y$ in $\G_2$  it is possible to reach in exactly two steps at most $23$ vertices of $\G_0$  in the same partite set as $y$, a contradiction as well (see Figure \ref{fig:G2G1G0} $(b)$).
\EndProof

From Proposition \ref{prop:G2G1G0} it immediately follows the main result of this section.

\begin{theorem}
\label{theo:bipartite(7,3,-4)}
There is no bipartite $(7,3,-4)$-graph.
\end{theorem}

Theorem \ref{theo:bipartite(7,3,-4)} settles the optimality of the known bipartite $(7,3,-6)$-graph, and therefore $\N^b(7,3)=80$.

%										Nb(11,3)>=190

\section{Three largest known bipartite graphs of diameter $3$}

In this section we present three new largest known bipartite graphs of degree $11$, diameter $3$ and order $190$. This improves by $4$ vertices the former lower bound for $\N^b(11,3)$.

To obtain such graphs we were inspired by  Observation \ref{obs:G1}, which tells us about the overall structure of a --hypothetical -- bipartite $(d,3,-4)$-graph $\G$ in the particular case of $\G_1$ being a spanning subgraph of $\G$ with exactly one connected component $\Phi_m$.

\begin{corollary}
Let $\G$ be a bipartite $(d,3,-4)$-graph such that $\G_1$ has exactly one connected component $G=\Phi_{d^2-d-1}$ and $V(G)$ spans $\G$. If the vertex $x_0$ in $G$ has neighbors $y_0,y_1,y_{-1},y_{i_1},y_{i_2},\ldots,y_{i_{d-3}}$ in $G$ then $x_k$ has neighbors $y_k,y_{k+1},y_{k-1},y_{k+i_1},y_{k+i_2},\ldots,y_{k+i_{d-3}}$ for every $k$.
\end{corollary}

When $d=4$ or $d=5$ we have as examples the existing graphs depicted in Figures \ref{fig:Bipartite(4,3,-4)} $(b)$ and \ref{fig:Bipartite(5,3,-4)}. It is then natural to ask if similar graphs exist for greater values of $d$.

\begin{problem}
Is there a a bipartite $(d,3,-4)$-graph with $d\ge 5$ such that $\G_1$ has exactly one connected component $G=\Phi_{d^2-d-1}$ and $V(G)$ spans $\G$?
\end{problem}

By computer search we obtained that for small degrees ($d=6,7,8,9$) such graphs do not exist. This is a strong indication that for all $d\ge 6$ the answer to the above problem is no. Thus, we shift our interest to a more general problem. 

We first introduce an extension to the construction of a $\Phi_m$. Let $d\ge 4$ and $a_1, a_2,\ldots, a_{d-3}$ be such that $2\le a_j\le m-2$ and $a_j\neq a_k$ when $j\neq k$. Then $\Phi_m(a_1, a_2,\ldots, a_{d-3})$ denotes the graph with vertex set $V=\{x_0,x_1,\ldots,x_{m-1}\}\cup\{y_0,y_1,\ldots,y_{m-1}\}$ and edge set $E=\{x_i\sim y_i, x_i\sim y_{i+1}, x_i\sim y_{i-1}, x_i\sim y_{i+a_j}|0\le i\le m-1$, $1\le j\le d-3\}$. As before, we do addition modulo $m$ on the vertex subscripts. Note that $\Phi_m(a_1, a_2,\ldots, a_{d-3})$ is a bipartite vertex-transitive graph.

\begin{problem}
\label{prob:largestm}
Given a natural number $d\ge 6$, find the largest natural number $m(d)$ for which there exist natural numbers $a_1, a_2,\ldots, a_{d-3}$ $(2\le a_j\le m-2)$ such that the graph $\Phi_{m(d)}(a_1, a_2,\ldots, a_{d-3})$ has diameter $3$.
\end{problem}

If we take a vertex $x_0$ of a  $\Phi_{m(d)}(a_1, a_2,\ldots, a_{d-3})$ and assume that $x_0$ has neighbors $y_0,y_1,y_{-1},y_{a_1},y_{a_2},\ldots,y_{a_{d-3}}$ then $x_0$ can reach in at exactly two steps the -- not necessarily distinct -- vertices $x_0,x_1,x_{-1},x_{2},x_{-2},x_{a_i},x_{-a_i}$, $x_{a_i+1}$, $x_{-a_i-1}$, $x_{a_i-1}$, $x_{-a_i+1}$ and $x_{a_i-a_j}$, and no other vertex. Since $\Phi_{m(d)}(a_1, a_2,\ldots, a_{d-3})$ is vertex-transitive, Problem \ref{prob:largestm} amounts to the following congruence-related problem:

\begin{problem}
Given a natural number $d\ge 6$, find the largest natural number $m(d)$ for which there exist natural numbers $a_1, a_2,\ldots, a_{d-3}$ such that the collection $0,1,-1,2,-2,a_i,-a_i,a_i+1,-a_i-1,a_i-1,-a_i+1,a_i-a_j$ of (not necessarily distinct) numbers  contains a full set of residues modulo $m(d)$.
\end{problem}

It is not difficult to verify that $m(d) \le d^2-d-1 = (\M^b(d,3)-4) / 2$.
    
With the aid of computer search we found the non-isomorphic bipartite $(11,3,-32)$-graphs $\Phi_{95}(4,7,16,27,38,52,62,81)$, $\Phi_{95}(4,16,30,43,51,62,71,89)$ and $\Phi_{95}(11,15,21,28,37,40,45,63)$. This discovery implies that $m(11)\ge 95$ and $\N^b(11,3)\ge 190$. Adjacency lists of these graphs are available at \url{http://guillermo.com.au/wiki/List_of_Publications} under the name of this paper.

%							CONCLUSIONS
\section{Conclusions}

In this paper we offered several structural properties for bipartite graphs of diameter $3$ and defect $4$. Using these properties we showed the non-existence of bipartite $(7,3,-4)$-graphs, which proves the optimality of the known bipartite $(7,3,-6)$-graph on $80$ vertices. This is just the second non-Moore bipartite graph known to be optimal.

We would also like to emphasize that, using the results of Section \ref{sec:Preliminaries} and reasoning as in Section \ref{sec:Nb(7,3)=80}, it is possible to prove as well the uniqueness of the only known bipartite $(5,3,-4)$-graph depicted in Figure  \ref{fig:Bipartite(5,3,-4)}, and the non-existence of bipartite $(6,3,-4)$-graphs.

In addition, some of the results in Section \ref{sec:Nb(7,3)=80} could have been stated for any bipartite $(d,3,-4)$-graph by providing a more elaborate proof. However, we decided to omit this extension as it does not lead to any conclusive outcome on the existence or otherwise of bipartite graphs of diameter $3$ and defect $4$ in general. We nevertheless feel that the following conjecture is valid.

\begin{conjecture}
There is no bipartite $(d,3,-4)$-graph with $d\ge 6$.
\end{conjecture}

%							BIBLIOGRAPHY
\def\cprime{$'$}
\providecommand{\bysame}{\leavevmode\hbox to3em{\hrulefill}\thinspace}
\providecommand{\MR}{\relax\ifhmode\unskip\space\fi MR }
\providecommand{\MRhref}[2]{%
  \href{http://www.ams.org/mathscinet-getitem?mr=#1}{#2}
}
\providecommand{\href}[2]{#2}


\begin{thebibliography}{10}


\bibitem{DJMP2}
C.~Delorme, L.~K. J{\o}rgensen, M.~Miller, and G.Pineda-Villavicencio, \emph{On
  bipartite graphs of diameter 3 and defect 2}, Journal of Graph Theory
  \textbf{61} (2009), no.~4, 271--288,
  \href{http://dx.doi.org/10.1002/jgt.20378}{doi:10.1002/jgt.20378}.

\bibitem{Die05}
R.~Diestel, \emph{{Graph Theory}}, 3rd. ed., Graduate Texts in Mathematics,
  vol. 173, Springer-Verlag, Berlin, 2005.

\bibitem{FH64}
W.~Feit and G.~Higman, \emph{The nonexistence of certain generalized polygons},
  Journal of Algebra \textbf{1} (1964), 114--131,
  \href{http://dx.doi.org/10.1016/0021-8693(64)90028-6}{doi:10.1016/0021-8693(%
64)90028-6}.

\bibitem{FPV12}
R.~Feria-Pur\'on and G.~Pineda-Villavicencio, \emph{On bipartite graphs of defect at most 4},
  Discrete Applied Mathematics \textbf{160} (2012), 140--154,
  \href{http://dx.doi.org/10.1016/j.dam.2011.04.018}{doi:10.1016/j.dam.2011.09.002}.

\bibitem{Jor93}
L.~K. J{\o}rgensen, \emph{Nonexistence of certain cubic graphs with small
  diameters}, Discrete Mathematics \textbf{114} (1993), 265--273,
  \href{http://dx.doi.org/10.1016/0012-365X(93)90371-Y}{doi:10.1016/0012-365X(%
93)90371-Y}.

\bibitem{LTS89}
C.~W.~H.~Lam, L.~Thiel, and S.~Swiercz, \emph{The nonexistence of finite projective planes of
order 10}, Canadian Journal of Mathematics {\bf 41} (1989), 1117--1123.

\bibitem{LH12} E.~Loz and P.~Hafner, {\it A bipartite $(7,3,-6)$-graph}, Personal communication with G.~Pineda-Villavicencio, 2010, \url{http://www.math.auckland.ac.nz/~hafner/bipartite/7.3}.

\bibitem{MS05a}
M.~Miller and J.~\v{S}ir\'{a}\v{n}, \emph{Moore graphs and beyond: {A} survey
  of the degree/diameter problem}, The Electronic Journal of Combinatorics
  \textbf{DS14} (2005), 1--61, dynamic survey.
  
\bibitem{Rys82}
H.~J.~Ryser, \emph{The existence of symmetric block designs}, Journal of Combinatorial Theory, Series A {\bf 32} (1982), \href{dx.doi.org/10.1016/0097-3165(82)90068-1}{doi:10.1016/0097-3165(82)90068-1}.


\end{thebibliography}
\end{document}